\documentclass[11pt,a4paper,reqno]{amsart}
\usepackage{enumitem,xcolor,graphicx,hyperref}\usepackage[cp1250]{inputenc}
\usepackage{ifthen,amssymb,amsmath,subfig}\usepackage{amsthm}\usepackage{amsfonts}
\newtheorem{theorem}{Theorem}[section]

\newtheorem{definition}{Definition}[section]

\newtheorem{corollary}{Corollary}[section]
\newtheorem{lemma}{Lemma}[section]
\newtheorem{remark}{Remark}[section]
\def\int{{\rm int}}\newcommand{\ds}{\displaystyle}
\allowdisplaybreaks

\begin{document}

\title[Second Hankel determinant for a certain subclass  ...]{Second Hankel determinant for a certain subclass of bi-close to convex functions defined by Kaplan}

\author[S. Kanas,  V. Sivasankari, O.Karthiyayini and S. Sivasubramanian]{S. Kanas$^{1}$, V.Sivasankari$^{2}$, O.Karthiyayini$^{2}$, and S. Sivasubramanian$^{3}$}

\address{$^{1}$University of Rzeszow, Al. Rejtana 16c, PL-35-959 Rzesz\'{o}w, Poland}\email{skanas@ur.edu.pl}
\address{$^{2}$Department of Science and Humanities, PES University-Electronic city Campus, Bangalore - 560100, Karnataka, India} \email{sivasankariv@pes.edu;  karthiyayini.roy@pes.edu}
\address{$^{3}$Department of Mathematics, University College of Engineering Tindivanam, Anna University (Chennai), Tindivanam 604001, Tamil Nadu, India} \email{sivasaisastha@rediffmail.com}

\thanks{$^{1}$Corresponding author}
\maketitle

\begin{abstract}
In this paper, we consider the class of strongly bi-close-to-convex functions of order $\alpha$ and  bi-close-to-convex functions of order $\beta$. We obtain an upper bound estimate for the second Hankel determinant for functions belonging to these classes. The results in this article improve some earlier result obtained  for the class of bi-convex functions.
\end{abstract}

\section{Introduction}\setcounter{equation}{0}

\subsection{Bi-univalence} If $f$ is in the class $\mathcal{S}$, then $f$ is one-to-one in $\mathbb{D}=\{z\in \mathbb{C}:\ |z|<1\}$, and
\begin{equation}\label{eq1.1} f(z)= z +  a_2z^2+a_3 z^3 +\cdots,
 \end{equation}
then the inverse $f^{-1}$ of $f$ has Maclaurin expansion in a disk of radius at least $1/4$, say
 \begin{equation}\label{eq1.2}
 f^{-1}(w)= w - a_{2}w^{2}+(2a_{2}^{2}-a_{3})w^{3}-(5a_{2}^{3}-5a_{2}a_{3}+a_{4})w^{4}+\cdots.\\
 \end{equation}
 An analytic function $f$ of the form \eqref{eq1.1} is said to be \textit{bi-univalent} in $\mathbb {D}$ if both $f$ and $f^{-1}$ are univalent in $\mathbb D$, in the sense that $f^{-1}$ has an univalent analytic continuation to $\mathbb{D}$. Let $\Sigma$ denote the class of all bi-univalent functions in $ \mathbb{D}$, given by the Taylor-Maclaurin series expansion \eqref{eq1.1}. Family $\Sigma$ has been the focus of attention for more than fifty years.  In \cite{Lew1}, Lewin  established  that for $f \in \Sigma$, $|a_2|\leq1.51$. Later on,  Brannan and Clunie \cite{BraClun} hypothesized that $|a_2|\leq\sqrt{2}$, however, their hypothesis has not been proved.  One of the results which deserves more attention but somehow unnoticed  is that of  Netanyahu \cite{Netan} who obtained a sharp upper bound $|a_2|\le\dfrac43$ for a class $\Sigma_1\subset \Sigma$, consisting of the functions that are  bi-univalent and its range contain $\mathbb{D}$. However, the sharp lower bound of the second coefficient $|a_2|$ in the class $\Sigma$ is not known, as well as bounds for succesive coefficients $|a_n|\ (n>2)$.
 Some examples of bi-univalent functions  are:\  $\ds\frac{z}{1-z}, \ds \frac12\log\left(\frac{1+z}{1-z}\right)$ or $\ds -\log(1-z)$, however the familiar Koebe function, or\ $\ds \frac{z}{1-z^2}$, which are the members of $\mathcal{S}$, are not the elements of the class $\Sigma$.

\subsection{Subfamilies of $\mathcal{S}$ and related bi-univalent functions}  Let $0 \leq \beta < 1$. The subclasses of $\mathcal{S}$ consisting of \textit{starlike functions of order }$\beta$ (and \textit{convex functions of order} $\beta$, respectively) are denoted by $\mathcal{ST}(\beta)$ ($\mathcal{CV}(\beta)$, resp.), and are defined analytically
\begin{equation}\label{eq1.3}
\mathcal{ST}(\beta)= \left\{ f \in \mathcal{S}: \Re \left(zf'(z)/f(z)\right)>\beta, \quad z \in \mathbb{D}\right\},
\end{equation}
\begin{equation}\label{eq1.4}
\mathcal{CV}(\beta)= \left\{f\in \mathcal{S}: \Re \left( 1+ zf''(z)/f'(z)\right)> \beta, \quad z \in \mathbb{D}\right\}.
\end{equation}
A function $f$ of the form \eqref{eq1.1} is said \textit{close-to-convex} in $\mathbb{D}$, if and only if there exists  a function $\Phi \in\mathcal{CV} = \mathcal{CV}(0)$ such that $\ds \Re \left( \frac{f'(z)}{\Phi'(z)}\right) > 0, z \in \mathbb{D}$.
 The family of normalized close-to-convex functions was first introduced  by  Kaplan \cite{ka52} and denoted $\mathcal{K}$.

Brannan and Taha \cite{bt86} introduced the classes $\mathcal{ST}_{\Sigma}(\beta)$ (and $\mathcal{CV}_{\Sigma}(\beta)$) of \textit{bi-starlike functions of order} $\beta$ (and \textit{bi-convex functions of order} $\beta$, resp.) corresponding to $\mathcal{ST}(\beta)$ and $\mathcal{CV}(\beta)$ defined by (\ref{eq1.3}) and (\ref{eq1.4}). They also found non-sharp estimates on $|a_{2}|$ and $|a_{3}|$ for its members of the form \eqref{eq1.1}. Following Brannan and Taha  \cite{bt86}, many researchers  (see, for example \cite{al12, hw12}) have recently introduced and investigated several interesting subclasses of $\Sigma$ and found non-sharp estimates on the first two Taylor-Maclaurin coefficients. Also, in \cite{bt86} the class of \textit{strongly bi-starlike functions of order} $\alpha$, where $0 < \alpha \leq 1$ has been defined and denoted  $\mathcal{ST}_{\Sigma}[\alpha]$. A function $f$ is in the class $\mathcal{ST}_{\Sigma}[\alpha]$,  if
$$\ds \left|\arg\left(\frac{zf^{\prime}(z)}{f(z)}\right)\right|<\frac{\alpha\pi}{2} \quad
\quad\textit{and}\quad
 \left|\arg\left(\frac{w g^{\prime}(w)}{g(w)}\right)\right|<\frac{\alpha\pi}{2} \quad (z, w \in \mathbb{D}),$$
where $g $ is the analytic continuation of $f^{-1}$ to $\mathbb{D}$.

For $0\leq \alpha \leq 1$, let $\mathcal{K}_{\alpha}$ denote the family of analytic and locally univalent In $\mathbb{D}$ functions $f$ of the form \eqref{eq1.1}, for which  there exists a convex function $\phi$ such that
\begin{equation}\label{eq 1.7}
\left|\arg \left(\frac{f^{\prime}(z)}{\phi^{\prime}(z)}\right)\right|<\frac{\alpha \pi}{2}\quad (z\in \mathbb{D}).
\end{equation}
That class has been introduced by Kaplan \cite{ka52} and later studied by Reade \cite{re56}. In particular, $\mathcal{K}_{0}$ is the family of convex univalent functions and $\mathcal{K}_{1}$ is the family of close-to-convex functions. Moreover, $\mathcal{K}_{\alpha_{1}}$ is a proper subclass of $\mathcal{K}_{\alpha_{ 2}}$ whenever $\alpha_{1}< \alpha_{2}$. An extension of $\mathcal{K}_{\alpha}$ is a  class $\mathcal{K}_{\alpha}(\beta)$ of close-to-convex functions of order $\beta$ \cite{re56}, given by
\begin{equation}\label{eq1.8}
\Re \left(\frac{f^{\prime}(z)}{\phi^{\prime}(z)}\right)>\beta \quad (z\in \mathbb{D}).
\end{equation}

Following Brannan and Taha the related families of bi-univalent functions have been considered, for example a class  $\mathcal{K}_{\Sigma}$ of bi-close-to-convex functions  \cite{si15}, a class of strongly bi-close-to-convex functions of order $\alpha$, denoted by    $\mathcal{K}_{\Sigma}[\alpha]$ and the class  of bi-close-to-convex functions of order $\beta$, denoted by $\mathcal{K}_{\Sigma}(\beta)$.

\begin{definition}\cite{si15}\label{def1.1}
 A function $f\in \Sigma$ of the form (\ref{eq1.1}) belongs to the class of bi-close to convex functions  $\mathcal{K}_{\Sigma}$, if there exist a function $\phi$, convex and univalent for $z \in \mathbb{D}$, such that
$$\Re \left\{\frac{f^{\prime}(z)}{\phi^{\prime}(z)}\right\} \geq 0, \quad \textit{and}\quad \Re \left\{\frac{g^{\prime}(w)}{\phi^{\prime}(w)}\right\} \geq 0 \quad  (z,w\in \mathbb{D}),$$
where $g$ is the analytic continuation of $f^{-1}$ to $\mathbb {D}$ with a series expansion  (\ref{eq1.2}).
\end{definition}
    \begin{definition}\label{de1.2}\cite{si15}
  Let $0 \leq \alpha \leq 1$. A function $f \in \Sigma$, given by (\ref{eq1.1}), is said to be strongly bi-close-to-convex of order $\alpha$ if there exist bi-convex functions $\phi$ and $\psi$ such that
 \begin{equation}\label{eq1.9}
 \left|\arg\left(\frac{f^{\prime}(z)}{\phi^{\prime}(z)}\right)\right|< \alpha \pi/2 \quad\textit{and}\quad
 \left|\arg\left (\frac{g^{\prime}(w)}{\psi^{\prime}(w)}\right)\right|< \alpha \pi/2 \quad (z, w \in \mathbb{D}).
    \end{equation}
Here, $g$ is the analytic continuation of $f^{-1}$ to $\mathbb{D}$. We denote the class of strongly bi-close-to-convex functions of order $\alpha$ by $K_{\Sigma}[\alpha]$.
\end{definition}

\begin{remark}\label{rem1.1}  We note that $\mathcal{K}_{\Sigma}[1] \equiv \mathcal{K}_{\Sigma}$, and  $\mathcal{K}_{\Sigma}[0] \equiv \mathcal{CV}_{\Sigma}$ \cite{bt86}.
\end{remark}

\begin{definition} \cite{si15} Let $0 \leq \beta < 1$. A function $f \in \Sigma$, given by (\ref{eq1.1}) is said to be  bi-close-to-convex of order $\beta $ if there exist the bi-convex functions $\phi$ and $\psi$  $\in \mathcal{CV}
 _{\Sigma}$ such that:
 \begin{equation}\label{eq1.14}
 \Re \left(\frac{f^{\prime}(z)}{\phi^{\prime}(z)}\right)>  \beta  \quad \textit{and}\quad
    \Re \left(\frac{g^{\prime}(w)}{\psi^{\prime}(w)}\right) > \beta  \quad (z, w \in \mathbb{D}),\end{equation}
where  $g$ is the analytic continuation of $f^{-1}$ to $\mathbb{D}$. We denote the class of  bi-close-to-convex functions of order $\beta$ by $\mathcal{K}_{\Sigma}(\beta)$.\end{definition}

\begin{remark}\label{rem1.2} We  note that $\mathcal{K}_{\Sigma}(0) \equiv \mathcal{K}_{\Sigma}$. Also, for  $\phi(z)=z$,  the class $N_{\Sigma}(\alpha)\quad(0 \leq \alpha \leq 1)$ reduces to the family of functions $f\in \Sigma$, satisfying the condition,
\begin{equation*}
 \left|\arg\,f^{\prime}(z)\right|< \alpha \pi/2 \quad \textit{and}\quad \left|\arg\,g^{\prime}(w)\right|< \alpha \pi/2 \quad (z, w \in \mathbb{D}),
\end{equation*}
and $\mathcal{K}_{\Sigma}(\beta)$ reduces to  $N_{\Sigma}(\beta)$ defined by the conditions
\begin{equation*}
 \Re \left(f^{\prime}(z)\right) >  \beta   \quad \textit{and}\quad  \Re \left(g^{\prime}(w)\right) > \beta  \quad (z, w \in \mathbb{D}),
\end{equation*}
where the function $g$ is defined by (\ref{eq1.2}). These classes were studied by \c{C}a\u{g}lar et al. \cite{Cag}
\end{remark}

Observe that if $f$ is given by (\ref{eq1.1}), then $g=f^{-1}$ is given by \eqref{eq1.2},  and if
 \begin{equation}\label{eq1.12}
 \phi(z)= z+c_{2}z^2+c_{3}z^3+c_{4}z^4+\cdots,
 \end{equation}
 then
 \begin{equation}\label{eq1.13}
 \psi(w)=\phi^{-1}(w)= w - c_{2}w^{2}+(2c_{2}^{2} - c_{3})w^{3}-(5c_{2}^{3}-5c_{2}c_{3}+c_{4})w^{4}+\cdots .
 \end{equation}
In the sequel we assume that $g,\phi,\psi$ have Taylor expansions as in \eqref{eq1.2}, \eqref{eq1.12} and \eqref{eq1.13}.

\subsection{Hankel determinant}  Towards the full understanding of a behavior of bi-univalence, it is necessary to extend our attention to  the Hankel determinants, that is one of the most important tool in Geometric Function Theory,  defined by Pommerenke \cite{Pom1, Pom2}. Noonan and Thomas \cite{jw76} defined the $q^{th}$ Hankel determinant of $f$ given by (\ref{eq1.1}) for integral $n\geq 1$ and $q \geq 1$ by

$$H_{q}(n)=\begin{vmatrix}
a_{n} & a_{n+1} & \cdots & a_{n+q-1}\\
a_{n+1} & a_{n+2} & \cdots & a_{n+q}\\
: & : & : & :\\
: & : & : & :\\
a_{n+q-1} & a_{n+q} & \cdots & a_{n+2q-2}
\end{vmatrix}$$\medskip

The importance of the Hankel determinants was recognized over half a century ago and it has been studied in great details, see for example \cite{Pom1, Pom2}.  The significance of the Hankel determinants follows from the study of singularities of analytic functions \cite[p. 329]{pd57}, see also \cite{Ed}, and  from the fact that it contains  the Fekete-Szeg\"o functional with its generalization \cite{Fek}.  Moreover,  $H_{2}(2)= a_{2}a_{4}-a_{3}^2$ is well known  second Hankel determinant. The Hankel determinant  is useful for estimating the modulus of coefficients and  the rate of growth of the coefficients. Both estimates determine the behavior of the studied function when the function itself and its properties are unknown.
Extensive studies of the Hankel determinant in the theory of meromorphic functions are due to Wilson \cite{Wil}; numerous applications in mathematical physics are given by Vein and Dale \cite{Vei}. Recently, many authors have discussed upper bounds for the Hankel determinant and Fekete-Szeg\"o functional for numerous subclasses of univalent functions \cite{pd57, Kan, fr69, jw76,  kl13, mot} and references therein. Very recently, the upper bounds of $H_{2}(2)$ for the classes $\mathcal{S}_{\Sigma}^*(\alpha)$ and $K_{\Sigma}(\alpha)$ were investigated by Deniz et al. \cite{ed15}, and extended by Orhan et al. \cite{ho16, ho18}.

Sivasubramanian et al. \cite{si15} found the  estimates for $|a_{2}|$ and $|a_{3}|$ for  the classes $\mathcal{K}_{\Sigma},\;\mathcal{K}_{\Sigma}[\alpha]$ and $ \mathcal{K}_{\Sigma}(\beta)$. Further they  verified Brannan and Clunie's conjecture $|a_{2}|\leq \sqrt{2}$ for some of their subclasses.

Therefore, a naturally arising problem addressed in this paper is to investigate the behavior of the Hankel determinants in the newly defined families.

\subsection{Some useful bounds} Let $\mathcal{P} $ denote the class of functions $p(z)$ of the form
  \begin{equation}\label{eq1.16}
  p(z)=1+p_{1}z+p_{2}z^2+p_{3}z^3+... ,
  \end{equation}
  which are analytic in the open unit disk $\mathbb{D}$ and such that $\Re\, p(z) >0,\ z\in \mathbb{D}$.

\begin{lemma}\label{1.4}\cite{cp75}
If the function  $ p \in \mathcal{P} $ is given by the series (\ref{eq1.16}), then  $|p_{k}| \leq 2,\ k=1,2,...$ .\end{lemma}

\begin{lemma}\label{1.5}\cite{ug58}
If the function $p \in \mathcal{P} $ is given by the series (\ref{eq1.16}), then\\
$$\begin{array}{rcl}
2 p_{2}&=& p_{1}^2 + x(4-p_{1}^2),\\
4 p_{3}&=& p_{1}^3 + 2p_{1}(4-p_{1}^2)x -p_{1}(4- p_{1}^2)x^2 + 2(4-p_{1}^2)(1-|x|^2)z,\end{array}$$
for some x, z with $|x| \leq 1 $ and $ |z| \leq 1$.\end{lemma}

 \begin{lemma}\label{1.6}\cite{ka12}
If the function $\phi \in \mathcal{CV}$, then for $\lambda \in \mathbb{R}$ ,
\begin{equation}\label{1.18}
|c_{3}-\lambda c_{2}^{2}| \leq \left \{ \begin{array}{rcl} 1-\lambda &\mbox{for} &  \lambda < 2/3, \\  1 &\mbox{for} &  2/3 \leq \lambda\leq 4/3,\\
\lambda -1 & \mbox{for}  & \lambda > 4/3. \end{array}\right.
\end{equation}
\end{lemma}

\begin{lemma}\label{1.7}\cite{ja07}
If the function $\phi \in \mathcal{CV}$, then $|c_{2}c_{4}-c_{3}^2| \leq  \frac{1}{8}.$\end{lemma}

\begin{lemma}\label{1.8}\cite{ba10} If the function $\phi \in \mathcal{CV}$, then $|c_{2}c_{3}-c_{4}| \leq \frac{1}{6}.$\end{lemma}

\section{Second Hankel determinant in class $\mathcal{K}_{\Sigma}[\alpha]$ and $\mathcal{K}_{\Sigma}(\beta)$} \label{sec2}\setcounter{equation}{0}
The first  aim of this section  is to find the best bound of the second Hankel determinant in the class $\mathcal{K}_{\Sigma}[\alpha]$.  A successful method of finding such bound has been exploited in \cite{si15} and other related publications.
\subsection{The class $\mathcal{K}_{\Sigma}(\beta)$} In the family of strongly bi-close-to-convex of order $\alpha$  we have the following non-sharp estimates of $H_2(2)$ however, this bound, for a particular selection of $\alpha$, improves the earlier results in \cite{ed15}.

\begin{theorem}\label{thm2.1}
Let $0 \leq \alpha \leq 1 $, and let the function $f$, given by \eqref{eq1.1}, be in the class $ \mathcal{K}_{\Sigma}[\alpha]$. Then
\begin{equation}\label{eq2.1}
|a_{2}a_{4}-a_{3}^2| \leq \frac{1}{8}+\frac{3}{2}\alpha +\frac{43}{9}\alpha^2 +\frac{1}{3}\alpha^3 +\frac{4}{3} \alpha^4.
\end{equation}\end{theorem}
\begin{proof}
From the condition \eqref{eq1.9} it follows that there exists $p, q \in \mathcal{P}$ such that
\begin{equation}\label{eq2.2}
f'(z)=\phi'(z)[p(z)]^{\alpha} \quad\textit{and}\quad g'(w)= \psi'(w)[q(w)]^{\alpha}.
\end{equation}
Let $p$ be given by \eqref{eq1.16} and  $q $ has a series representations
\begin{equation}\label{eq2.5}
q(w)= 1 + q_{1}w+ q_{2}w^2 + q_{3}w^3+\cdots.
\end{equation}
Then, equating the coefficients of both sides of \eqref{eq2.2}, when $f, p, q, \phi$ and $\psi$ have given power series, we obtain a number of equalities, below.
\begin{equation}\label{eq2.6}
2a_{2}= 2c_{2}+\alpha p_{1},\quad -2a_{2}= -2c_{2}+\alpha q_{1},
\end{equation}
\begin{equation}\label{eq2.7}
3a_{3}= 3c_{3}+ 2\alpha c_{2} p_{1} + \alpha p_{2}+ \frac{1}{2} \alpha (\alpha-1)p_{1}^2,
\end{equation}
\begin{equation}\label{eq2.8}
4a_{4} = 4c_{4}+3c_{3} \alpha p_{1}+2 \alpha c_{2}p_{2} +\alpha p_{3}+\alpha(\alpha-1)c_{2}p_{1}^2 + \alpha(\alpha-1)p_{1}p_{2} + \frac{\alpha(\alpha-1)(\alpha-2)}{6}p_{1}^3,
\end{equation}
\begin{equation}\label{eq2.10}
6a_{2}^2 - 3a_{3}= 6c_{2}^2 - 3c_{3}- 2\alpha c_{2} q_{1} + \alpha q_{2}+ \frac{1}{2} \alpha (\alpha-1)q_{1}^2,
\end{equation}
\begin{equation}\label{eq2.11}
\begin{split}
& -20 a_{2}^3+ 20 a_{2}a_{3} - 4 a_{4} = -20 c_{2}^3 +20 c_{2}c_{3}
 -4c_{4}+ 6c_{2}^2 \alpha q_{1}- 3c_{3} \alpha q_{1}- 2 \alpha c_{2}q_{2}\\
& \quad \quad \quad \quad \quad \quad \quad  +\alpha q_{3}- \alpha(\alpha-1)c_{2}q_{1}^2 + \alpha(\alpha-1)q_{1}q_{2}+ \frac{\alpha(\alpha-1)(\alpha-2)}{6}q_{1}^3 .\\
\end{split}
\end{equation}

The equality  \eqref{eq2.6} immediately gives $p_{1}= - q_{1}$. Next, by (\ref{eq2.7}) and (\ref{eq2.10}), we obtain
\begin{equation}\label{eq2.12}
a_{3}= c_{3}+ \alpha c_{2}p_{1} + \frac{1}{4} \alpha^2 p_{1}^2  + \frac{\alpha}{6} ( p_{2}-q_{2}).
\end{equation}
Similarly, making necessary calculations of (\ref{eq2.8}) and (\ref{eq2.11}), we get
\begin{equation}\label{eq2.13}
\begin{split}
& a_{4}= c_{4} + \frac{\alpha}{8}(p_{3}-q_{3})+\frac{1}{4}
c_{2} \alpha (p_{2}+q_{2}) +\frac{5}{24} \alpha^2 p_{1}(p_{2}-q_{2})\\
& \quad \quad \quad \quad  +\frac{5}{12}\alpha c_{2}(p_{2}-q_{2})+ \frac{5}{4} c_{3} \alpha p_{1}+ \frac{1}{4} \alpha(\alpha-1) c_{2}p_{1}^2\\
& \quad \quad \quad \quad +\frac{1}{8} \alpha (\alpha -1)p_{1}(p_{2}+q_{2})-\frac{1}{2}c_{2}^2\alpha p_{1}+\frac{1}{24}\alpha (\alpha-1)(\alpha-2) p_{1}^3.
\end{split}
\end{equation}
Hence
\begin{equation}\label{eq2.14}
\begin{split}
& \left|a_{2}a_{4}-a_{3}^2\right| = \left|c_{2}c_{4}- c_{3}^2 +\frac{1}{2} \alpha c_{4}p_{1}
-\frac{3}{4}c_{2}c_{3}\alpha p_{1} +\frac{1}{8}c_{3} \alpha^{2} p_{1}^2- \frac{5}{4}\alpha^{2} p_{1}^2 c_{2}^2 \right.\\
& \quad \quad \quad \quad   -\frac{1}{3} c_{3}\alpha (p_{2}-q_{2})+ \frac{5}{12}c_{2}^{2} \alpha (p_{2}- q_{2})- \frac{1}{2}c_{2}^3 \alpha p_{1} -\frac{1}{16}\alpha^4 p_{1}^4\\
& \quad \quad \quad \quad  + \frac{1}{12}c_{2} \alpha^2 p_{1}(p_{2}- q_{2})
 +\frac{1}{48}\alpha^2(\alpha-1)(\alpha-2)p_{1}^4\\
 & \quad \quad \quad \quad +\frac{1}{8}\alpha^2 (\alpha-1)c_{2}p_{1}^3+ \frac{1}{24} c_{2} \alpha(\alpha-1)(\alpha-2)p_{1}^3\\
&\quad \quad \quad \quad  + \frac{1}{8}\alpha c_{2}(p_{3}-q_{3})+\frac{1}{8}c_{2}\alpha^2 p_{1}(p_{2}+q_{2}) \\
&\quad \quad \quad \quad+\frac{1}{48}\alpha^3p_{1}^2(p_{2}-q_{2})+\frac{1}{4}\alpha(\alpha-1)c_{2}^2 p_{1}^2 -\frac{1}{2}c_{2}\alpha^3p_{1}^3\\
&\quad \quad \quad \quad +\frac{1}{16} \alpha^2 p_{1}(p_{3}-q_{3}) +\frac{1}{16} \alpha^2 (\alpha-1)p_{1}^2(p_{2}+q_{2})\\
&\quad \quad \quad \quad + \left.\frac{1}{4}c_{2}^2 \alpha(p_{2}+q_{2}) + \frac{1}{8} c_{2} \alpha(\alpha-1)  p_{1}(p_{2}+q_{2}) - \frac{1}{36} \alpha^{2} (p_{2}-q_{2})^{2}\right|.
\end{split}
\end{equation}
Let's apply Lemma \ref{1.5} to $p_2$ and $q_2$. Then, for some $x,y$ such that $|x|\le 1,\ |y|\le 1$, it holds
$$ 2p_{2}= p_{1}^2+x(4-p_{1})^2, \quad 2q_{2}= q_{1}^2 +y (4-q_{1})^2. $$
from which we have
\begin{equation}\label{eq2.15}p_{2}-q_{2}= \frac{(4-p_{1}^2)(x-y)}{2} ; \quad p_{2}+q_{2}= p_{1}^2+\frac{(4-p_{1}^2)(x+y)}{2}.
\end{equation}
Apply now Lemma \ref{1.5} to $p_3$, $q_3$ and  obtain
\begin{equation}\label{eq2.16}p_{3}-q_{3}= \frac{p_{1}^3}{2} + \frac{p_1(4-p_{1}^2)(x+y)}{2}- \frac{p_1(4-p_{1}^2)(x^2 +y^2)}{4}
 +\frac{(4-p_{1}^2)}{2}\left[(1-|x|^2)z-(1-|y|^2)w \right],
\end{equation}
for some $ x,y,z$ and $w $ with $ |x| \leq 1$, $|y| \leq 1,$ $|z|\leq 1 $ and $ |w| \leq 1.$
Making use of  \eqref{eq2.15} and \eqref{eq2.16} to \eqref{eq2.14} gives
 \begin{equation}\label{eq2.17}
\begin{split}
&\left|a_{2}a_{4}-a_{3}^2\right| = \left|c_{2}c_{4}- c_{3}^2+ \frac{\alpha}{2}(c_{4}-c_{2}c_{3})p_{1} -\frac{1}{4}\alpha c_{2}c_{3}p_{1}+ \frac{1}{8} \alpha^2 (c_{3}-10c_{2}^2)p_{1}^2 \right.\\
&\quad \quad \quad \quad -\frac{1}{3}\alpha \left(c_{3}- \frac{5}{4}c_{2}^2 \right)(p_{2}-q_{2})+ \frac{1}{12}c_{2} \alpha^2 p_{1}(p_{2}- q_{2})
- \frac{1}{2}c_{2}^3 \alpha p_{1} - \frac{1}{2}c_{2}\alpha^3p_{1}^3 \\ & \quad \quad \quad \quad -\frac{1}{16}\alpha^4 p_{1}^4 +\frac{1}{16} \alpha^2(\alpha -1) p_{1}^2 (p_{2}- q_{2})+\frac{1}{4}\alpha(\alpha-1)c_{2}^2 p_{1}^2
 +\frac{1}{48} \alpha^2(\alpha-1)(\alpha-2)p_{1}^4
\\ & \quad \quad \quad \quad  +\frac{1}{8}\alpha^2(\alpha-1)c_{2}p_{1}^3 +\frac{1}{24}c_{2}\alpha(\alpha-1)(\alpha-2)p_{1}^3
  + \frac{1}{8}\alpha c_{2}(p_{3}-q_{3}) \\
&\quad \quad \quad \quad +\frac{1}{8}c_{2}\alpha(\alpha-1) p_{1}(p_{2}+q_{2})+\frac{1}{8}c_{2}\alpha^2 p_{1}(p_{2}+q_{2})+\frac{1}{48}\alpha^3 p_{1}^2(p_{2}-q_{2})\\
&\quad \quad \quad \quad + \left. \frac{1}{16} \alpha^2 p_{1}(p_{3}-q_{3})
 + \frac{1}{4}c_{2}^2 \alpha(p_{2}+q_{2})  - \frac{1}{36} \alpha^{2} (p_{2}-q_{2})^{2}\right|.
\end{split}
\end{equation}
Without lost of generality, we can restrict our considerations to $p_1:=p \in [0,2]$. Applying this and the triangle inequality to \eqref{eq2.17}, we have
\begin{equation}\label{eq2.19}
\begin{split}
&\left|a_{2}a_{4}-a_{3}^2\right| \leq \left|c_{2}c_{4}- c_{3}^2\right| + \frac{\alpha}{2}|(c_{4}-c_{2}c_{3})|p +\frac{1}{2}|c_{2}^3|\alpha p+\frac{\alpha}{4} |c_{2}c_{3}|p\\
&\quad \quad\quad\quad  + \frac{1}{8}\alpha^2|c_{3}-10c_{2}^2|p^2 +\frac{1}{2}|c_{2}|\alpha^3 p^3+\frac{1}{16}\alpha^4 p^4+\frac{1}{4}|\alpha(\alpha-1)| |c_{2}^2| p^2 \\
&\quad \quad \quad \quad+\frac{1}{48}|\alpha^2(\alpha-1)(\alpha-2)|p^4+ \frac{1}{8}|\alpha^2(\alpha-1)| |c_{2}|p^3+ \frac{1}{24}|c_{2}||\alpha(\alpha-1)(\alpha-2)|p^3\\
&\quad \quad \quad \quad +\frac{1}{16}|\alpha^2(\alpha-1)|p^4 +\frac{1}{16}\alpha| c_{2}|p^3+\frac{1}{8}|c_{2}|\alpha^2 p^3+\frac{1}{8}|c_{2}||\alpha(\alpha-1)|p^3+\frac{1}{32}\alpha^2 p^4\\
&\quad \quad \quad \quad +\frac{1}{4}|c_{2}^2| \alpha p^2 +(|x|+|y|) \left[\frac{\alpha}{6}\left|c_3- \frac{5}{4}c_{2}^2\right|+\frac{1}{8}|c_{2}^2| \alpha \right](4-p^2)\\
& \quad \quad \quad \quad + \left[\frac{1}{16}|c_{2}|\alpha+\frac{5}{48}|c_{2}|\alpha^2  +\frac{1}{16}|c_{2}||\alpha(\alpha-1)| \right] p(4-p^2)\\
&\quad \quad \quad \quad + \left.\left[\frac{1}{32}\alpha ^2+  \frac{1}{96}\alpha^3+\frac{1}{32}|\alpha^2(\alpha-1)| \right]p^2(4-p^2)\right]\\
& \quad \quad \quad \quad +(|x|^2+|y|^2)\left[\frac{\alpha^2}{64}p^2(4-p^2) +\frac{\alpha}{32}| c_{2}|p(4-p^2)\right]\\
&\quad \quad \quad \quad +(1-|x|^2)\left[ \frac{\alpha}{16}|c_{2}|(4-p^2)+ \frac{\alpha^2}{32}p(4-p^2)\right]\\
&\quad \quad \quad \quad + (1-|y|^2)\left[ \frac{\alpha}{16}|c_{2}|(4-p^2)+ \frac{\alpha^2}{32}p(4-p^2)\right]\\
&\quad \quad \quad \quad + \frac{\alpha^2}{144}(4-p^2)^2(|x|+|y|)^2.
\end{split}
\end{equation}
\noindent We now apply Lemma \ref{1.6}, Lemma \ref{1.7} with Lemma \ref{1.8} to \eqref{eq2.19}, and  deduce that
{\allowdisplaybreaks
{\begin{eqnarray*}
\left|a_{2}a_{4}-a_{3}^2\right|& \leq & \frac{1}{8} + \frac{5}{6}\alpha p+
\left(\frac{9}{8}\alpha^2+  \frac{1}{4}\alpha(1-\alpha)+\frac{1}{4}\alpha\right)p^2 \\
 && + \left(\frac{1}{16}\alpha +\frac{1}{8}\alpha^2+\frac{1}{2}\alpha^3 +\frac{1}{8}\alpha(1-\alpha)\frac{2\alpha+5}{3} \right) p^3 \\
 && +  \left[\frac{1}{16} \alpha^4  +\frac{1}{32}\alpha^2 +\frac{\alpha^2(1-\alpha)(5-\alpha)}{48}\right]p^4\\
&& +  \frac{(\alpha)}{8}(4-p^2)+\frac{(\alpha)^2}{16}p(4-p^2)\\
 && +  (|x|+|y|)\left[ \left( \frac{\alpha^3}{96} +\frac{\alpha^2}{32}+\frac{1}{32}\alpha^2(1-\alpha)\right) p^2(4-p^2)\right.\\
 && +  \left.\left(\frac{5\alpha^2}{48}+\frac{\alpha}{16} +\frac{\alpha(1-\alpha)}{16}\right)p(4-p^2)+\frac{\alpha}{6}(4-p^2)\right]\\
&& + (|x|^2+|y|^2)\left[ \frac{\alpha^2}{64}(p^2-4)(4-p^2)\right]+\frac{\alpha^2}{144}(4-p^2)^2 (|x|+|y|)^2.
\end{eqnarray*}}}
Taking $ \gamma_{1} =|x|\leq 1,\gamma_{2}=|y| \leq 1,$ we  rewritte the above as follows
$$|a_{2}a_{4}-a_{3}^2| \leq S_{1} + S_{2}(\gamma_{1} +\gamma_{2})+S_{3}(\gamma_{1}^2+\gamma_{2}^2)+S_{4}(\gamma_{1}+\gamma_{2})^2= F(\gamma_{1},\gamma_{2}),$$ where
\begin{eqnarray*}
S_{1} = S_{1}(p) &= &\frac{1}{8} + \frac{5}{6}\alpha p+
\left(\frac{9}{8}\alpha^2+  \frac{1}{4}\alpha(1-\alpha)+\frac{1}{4}\alpha\right)p^2\\
 &+& \left(\frac{1}{16}\alpha +\frac{1}{8}\alpha^2+\frac{1}{2}\alpha^3 +\frac{1}{8}\alpha(1-\alpha)\frac{2\alpha+5}{3} \right) p^3 \\
 &+& \left[\frac{1}{16} \alpha^4  +\frac{1}{32}\alpha^2 +\frac{\alpha^2(1-\alpha)(5-\alpha)}{48}\right]p^4\\
& +& \frac{(\alpha)}{8}(4-p^2)+\frac{(\alpha)^2}{16}p(4-p^2) \geq 0\\
S_{2} = S_{2}(p) &= & \left( \frac{\alpha^3}{96} +\frac{\alpha^2}{32}+\frac{1}{32}\alpha^2(1-\alpha)\right) p^2(4-p^2)\\
  &+&\left(\frac{5\alpha^2}{48}+\frac{\alpha}{16} +\frac{\alpha(1-\alpha)}{16}\right)p(4-p^2)+\frac{\alpha}{6}(4-p^2) \geq 0 \\
S_{3}= S_{3}(p)& =&\displaystyle \frac{\alpha^2}{64}(p^2-4)(4-p^2) \leq 0 \\
S_{4}= S_{4}(p)&=& \displaystyle\frac{\alpha^2}{144}(4-p^2)^2 \geq 0 \\
\end{eqnarray*}
In order to obtain an estimate of $|H_2(2)|$  we need to maximize $F(\gamma_{1},\gamma_{2})$ in the closed square
$$\Delta:=\{(\gamma_{1},\gamma_{2}): 0 \leq \gamma_{1} \leq 1,0 \leq \gamma_{2} \leq 1\}.$$
Since  $S_{3}<0$ and $ S_{3}+ 2S_{4}>0$  and $p\in (0,2)$ , we conclude that $F_{\gamma_{1}\gamma_{1}}F_{\gamma_{2}\gamma_{2}}-(F_{\gamma_{1}\gamma_{2}})^2 < 0$\ for all $\gamma_{1},\gamma_{2} \in \int\Delta,$
and thus the function $F$ can attain  a maximum only  on the boundary of  $\Delta$.

We first note that $F$ is symmetric in $\gamma_1$ and $\gamma_2$, therefore it is enough to consider $0\le \gamma_1\le 1$ and $0\le \gamma_2\le \gamma_1$. For $ \gamma_{2}=0 $ and $0 \leq \gamma_{1} \leq 1$, we obtain
$$F(\gamma_{1},0)= G(\gamma_{1})= S_{1}+S_{2}\gamma_{1}+(S_{3}+S_{4})\gamma_{1}^2.$$

Fix $p\in [0,2]$ and  consider two separate cases:\\

\noindent \textbf{(i)} $S_{3}+S_{4} \geq 0$. In this case $ G'(\gamma_{1})= 2(S_{3}+S_{4})\gamma_{1}+S_{2}>0$,  that is, $G(\gamma_{1})$ is an increasing function. Hence for fixed $p \in [0,2)$  the maximum of $G(\gamma_{1})$ may occurs only at $\gamma_{1}=1 $, and
$$\max\, G(\gamma_{1})= G(1)= S_{1}+S_{2}+S_{3}+S_{4}.$$

\noindent \textbf{(ii)} $S_{3}+S_{4}<0$.  Since $S_{2}+2(S_{3}+S_{4}) \geq 0$ for $0 <\gamma_{1}<1$,  it is clear that   $S_{2}+2(S_{3}+S_{4})< 2(S_{3}+S_{4})\gamma_{1}+S_{2}<S_{2}$  so that $ G'(\gamma_{1})>0$. Hence, similarly as in the case (i)  the maximum of $G(\gamma_{1})$  is attained for $\gamma_{1}=1.$\\

For $\gamma_{1}=1$ and $0 \leq \gamma_{2} \leq 1$, we obtain
$$F(1,\gamma_{2})= H(\gamma_{2})= (S_{3}+S_{4})\gamma_{2}^2+(S_{2}+2S_{4})\gamma_{2} +S_{1}+S_{2}+S_{3}+S_{4}.$$
Similarly, to the above cases of $S_{3}+S_{4},$ we get that
$$\max\ H(\gamma_{2})=H(1)= S_{1}+2S_{2}+2S_{3}+4S_{4}.$$
Since $G(1)\leq H(1)$ for $p \in [0,2]$ we have that $\max\, F(\gamma_{1}, \gamma_{2})= F(1,1)$ on the boundary of $\Delta$ and from this on the closed square $\Delta$.

Next, let's define a function $ K: [0,2] \rightarrow \mathbb{R}$ as follows
\begin{equation}\label{eq2.20}
K(p)= \max\, F(\gamma_{1},\gamma_{2})= F(1,1)= S_{1}+2S_{2}+2S_{3}+4S_{4},
\end{equation}
that is, in view of \eqref{eq2.20},
\begin{eqnarray*}
K(p) &= &\frac{1}{8} + \frac{11}{6}\alpha -\frac{1}{18}\alpha^2+\left(\frac{11}{6}\alpha +\frac{7}{12}\alpha^2 \right)p
 \\
&  + &\left(\frac{1}{24}\alpha +\frac{101}{72}\alpha^2-\frac{1}{6}\alpha^3 \right) p^2 \\
& + &\left[\frac{5}{12} \alpha^3 -\frac{19}{48}\alpha^2 +\frac{25}{48}\alpha \right]p^3\\
& + &\left[\frac{1}{12} \alpha^4 -\frac{7}{48}\alpha^3 +\frac{11}{144}\alpha^2+\frac{1}{16}\alpha \right]p^4.
 \end{eqnarray*}

By an elementary calculation, we find that
\begin{eqnarray*}
K'(p) &= & \left(\frac{11}{6}\alpha +\frac{7}{12}\alpha^2 \right) + \left(\frac{1}{12}\alpha +\frac{101}{36}\alpha^2-\frac{1}{3}\alpha^3 \right) p \\
& + &\left[\frac{5}{4} \alpha^3 -\frac{19}{16}\alpha^2 +\frac{25}{16}\alpha \right]p^2
 + \left[\frac{1}{3} \alpha^4 -\frac{7}{12}\alpha^3 +\frac{11}{36}\alpha^2+\frac{1}{4}\alpha \right]p^3,
\end{eqnarray*}
that can be rewritten as
\begin{eqnarray*}
K'(p) &= & \left(\frac{11}{6}\alpha +\frac{7}{12}\alpha^2 \right) + \frac{\alpha\, p}{3}\left[\alpha(1-\alpha) +\frac{89}{22}\alpha+\frac{1}{4} \right] \\
& + &\frac{\alpha\, p^2}{4}\left[5\alpha^2 +\frac{19}{4}(1-\alpha) +\frac{3}{2} \right] + \frac{\alpha\, p^3}{3}\left[\alpha\left(\alpha^2-\alpha+\frac{11}{12}\right)+\frac{3}{4}(1-\alpha^2) \right],
\end{eqnarray*}
from which it is easily seen that  $K'(p)>0$ for $0 <  \alpha \leq 1$. Hence  $K(p)$  is an increasing function of $p$ so that $K(p)$  attains its maximum at $p=2$. Consequently, we have
$$ \max_{0 \leq p \leq 2}\, K(p)= K(2) = \frac{1}{8}+\frac{3}{2}\alpha +\frac{43}{9}\alpha^2+ \frac{1}{3} \alpha^3 +\frac{4}{3}\alpha^4.$$
This completes the proof of the theorem.
\end{proof}

\begin{remark} For $\alpha = 1$, we have the following bound
$$|a_2 a_4 - a_3^{2} |\leq \frac{581}{72},$$
and when $\phi(z)=z$, Theorem \ref{thm2.1} reduces to the Theorem 2 in \cite{Cag}. Also, when
$\alpha =0$, we get the estimate for the class of bi-convex functions, which significantly improves the bound due to Deniz et al. \cite{ed15}, below. Unfortunately, we do not know if that result is sharp.\end{remark}

\begin{corollary} For $0 \leq \alpha < 1 $, and $f\in \mathcal{K}_{\Sigma}\equiv \mathcal{K}_{\Sigma}[0]$ we have
\begin{equation}
|a_{2}a_{4}-a_{3}^2| \leq \frac{1}{8}.
\end{equation}
\end{corollary}

\subsection{The class $\mathcal{K}_{\Sigma}(\beta)$}\label{sec3} In order to estimate the second Hankel determinat in $\mathcal{K}_{\Sigma}(\beta)$ we apply consideration similar  to that used in the proof of \ref{thm2.1}.

\begin{theorem}\label{thm3.1}
Let $0 \leq \beta < 1 $, and let the function $f$ given by (\ref{eq1.1}) be in the class $ \mathcal{K}_{\Sigma}(\beta)$. Then
\begin{equation}\label{eq3.1}
|a_{2}a_{4}-a_{3}^2| \leq \frac{1}{8}+ \frac{19}{6}(1-\beta)+6(1-\beta)^2+ 4(1-\beta)^3+(1-\beta)^4.
\end{equation}
\end{theorem}

\begin{proof}
\noindent By  \eqref{eq1.14}  there exist $p, q \in \mathcal{P}$ such that
\begin{equation}\label{eq3.2}
\frac{f'(z)}{\phi'(z)}= \beta + (1-\beta)[p(z)]\quad and \quad \frac{g'(w)}{ \psi'(w)}= \beta + (1- \beta)[q(w)].
\end{equation}
Let $p,q$ have series representations as in the previous section.  Then, equating coefficients of $z, z^2$ and $z^3$ of both sides of \eqref{eq3.2}, we obtain
\begin{equation}\label{eq3.6}
2a_{2}= 2c_{2}+ (1-\beta) p_{1}\quad and \quad -2a_{2}= -2c_{2}+(1-\beta) q_{1},
\end{equation}
\begin{equation}\label{eq3.7}
3a_{3}= 3c_{3}+ 2(1-\beta) c_{2} p_{1} + (1-\beta)p_{2},
\end{equation}
\begin{equation}\label{eq3.8}
4a_{4} = 4c_{4}+3c_{3} (1-\beta)p_{1}+2 (1-\beta) c_{2}p_{2} + (1-\beta) p_{3}.
\end{equation}
\begin{equation}\label{eq3.10}
6a_{2}^2 - 3a_{3}= 6c_{2}^2 - 3c_{3}- 2 (1-\beta) c_{2} q_{1} + (1-\beta) q_{2},
\end{equation}
\begin{equation}\label{eq3.11}
\begin{split}
& -20 a_{2}^3+ 20 a_{2}a_{3} - 4 a_{4} = -20 c_{2}^3 +20 c_{2}c_{3}
 -4c_{4}+ 6c_{2}^2 (1-\beta) q_{1}- 3c_{3} (1-\beta) q_{1} \\
&\quad \quad \quad \quad \quad \quad \quad \quad \quad \quad \quad - 2 (1-\beta) c_{2}q_{2} +(1-\beta) q_{3}.
\end{split}
\end{equation}
From \eqref{eq3.6} we get $p_{1}= - q_{1}$, and $a_{2}= c_{2}+\frac{(1-\beta)p_{1}}{2}$, and making use of  (\ref{eq3.7}), (\ref{eq3.10}), and \eqref{eq3.11},  we  have
\begin{equation}\label{eq3.12}
a_{3}= c_{3} + (1-\beta) c_{2}p_{1} + \frac{(1-\beta)(p_{2}-q_{2})}{6}+ \frac{(1-\beta)^2 p_{1}^2}{4},
\end{equation}
\begin{equation}\label{eq3.13}
\begin{split}
& a_{4}= c_{4} + \frac{1}{8}(1-\beta)(p_{3}-q_{3}) + \frac{1}{4}c_{2}(1-\beta)(p_{2}+q_{2})- \frac{1}{2}c_{2}^2(1-\beta)p_{1}\\
&\quad \quad\quad\quad +\frac{5}{4}c_{3}(1-\beta)p_{1}+ \frac{5}{12}c_{2}(1-\beta)(p_{2}-q_{2}) +\frac{5}{24}(1-\beta)^2 p_{1}(p_{2}-q_{2}).
\end{split}
\end{equation}

\noindent Hence
\begin{equation}\label{eq3.14}
\begin{split}
&\left|a_{2}a_{4}- a_{3}^{2} \right| = \left|c_{2}c_{4}- c_{3}^2
+\frac{1}{2}(1-\beta)c_{4} p_{1}-\frac{3}{4}c_{2}c_{3}(1-\beta)p_{1}\right.\\
& \quad \quad \quad \quad +\frac{1}{8}c_{3}(1-\beta)^{2} p_{1}^2 - \frac{5}{4}(1-\beta)^{2} p_{1}^2 c_{2}^2 -\frac{1}{3} c_{3}(1-\beta)(p_{2}-q_{2})\\
& \quad \quad \quad \quad + \frac{5}{12}c_{2}^{2} (1-\beta)(p_{2}- q_{2})- \frac{1}{2}c_{2}(1-\beta)^{3} p_{1}^3 + \frac{1}{12}c_{2}(1-\beta)^2 p_{1}(p_{2}- q_{2})\\
& \quad \quad \quad \quad - \frac{1}{2}c_{2}^3 (1-\beta) p_{1} -\frac{1}{16}(1-\beta)^4 p_{1}^4 + \frac{1}{48}(1-\beta)^3 p_{1}^2 (p_{2}- q_{2})\\
&\quad \quad \quad \quad + \frac{1}{16} (1-\beta)^{2} p_{1}(p_{3}- q_{3}) + \frac{1}{8}(1-\beta)c_{2}(p_{3}-q_{3})\\
&\quad \quad \quad \quad + \left. \frac{1}{4}c_{2}^2 (1-\beta)(p_{2}+q_{2}) + \frac{1}{8} c_{2} (1-\beta)^{2}  p_{1}(p_{2}+q_{2}) - \frac{1}{36} (1-\beta)^{2} (p_{2}-q_{2})^{2}\right|.
\end{split}
\end{equation}
Now, we apply  the relations  \eqref{eq2.15} and \eqref{eq2.16}  to \eqref{eq3.14} and we obtain
{\allowdisplaybreaks
{ \begin{eqnarray}\label{eq3.18}
\left|a_{2}a_{4}-a_{3}^2 \right| &= &\left| c_{2}c_{4}- c_{3}^2 + \frac{(1-\beta)}{2}(c_{4}-c_{2}c_{3})p_{1} -\frac{1}{4}(1-\beta) c_{2}c_{3}p_{1} \right. \nonumber \\
&& +  \frac{1}{8}(1-\beta)^2(c_{3}-10c_{2}^2) p_{1}^2 -\frac{1}{3}(1-\beta)\left(c_{3}- \frac{5}{4}c_{2}^2 \right)\left(\frac{4-p_{1}^2}{2}(x-y)\right) \nonumber \\
&& -  \frac{1}{2}c_{2}(1-\beta)^3 p_{1}^3+\frac{1}{12}c_{2}(1-\beta)^2 p_{1} \left(\frac{4-p_{1}^2}{2}(x-y)\right) \nonumber \\
&& - \frac{1}{2}c_{2}^3 (1-\beta) p_{1} -\frac{1}{16}(1-\beta)^4 p_{1}^4 + \frac{1}{48}(1-\beta)^3 p_{1}^2 \left(\frac{4-p_{1}^2}{2}(x-y)\right)\nonumber \\
&& +  \frac{1}{16} (1-\beta)^2 p_{1} \left[\frac{p_{1}^3}{2} +\frac{(4-p_{1}^2)p_{1}}{2}(x+y) \right. \nonumber \\
&& - \left. \frac{(4-p_{1}^2)p_{1}}{4}(\textcolor{white}{l}x^2+y^2\textcolor{white}{l})+ \frac{4-p_{1}^2}{\textcolor{white}{l}2\textcolor{white}{l}}\left[\textcolor{white}{l}(1-|x|^2)z-(1-|y|^2)w\textcolor{white}{l} \right]\right] \nonumber \\
&& +  \frac{1}{8}(1-\beta)c_{2} \left[\frac{p_{1}^3}{2} +\frac{(4-p_{1}^2)p_{1}}{2}(x+y) \right. \nonumber \\
&& -  \left. \frac{(4-p_{1}^2)p_{1}}{4}(x^2+y^2) + \frac{4-p_{1}^2}{2}\left[(1-|x|^2)z-(1-|y|^2)w \right] \right] \nonumber \\
&& + \frac{1}{4}c_{2}^2 (1-\beta)\left(p_{1}^2+\frac{4-p_{1}^2}{2}(x+y)\right) + \frac{1}{8} c_{2} (1-\beta)^2 p_{1}\left(p_{1}^2+\frac{4-p_{1}^2}{2}(x+y)\right) \nonumber \\
&& - \left.\frac{1}{36} (1-\beta)^2 \frac{\textcolor{white}{l}(4-p_{1}^2)^2 (x-y)^2\textcolor{white}{l}}{\textcolor{white}{l}4\textcolor{white}{l}} \right|,
\end{eqnarray}}}
where $x,y,z$ and $w$ are such that $ |x| \leq 1$, $|y| \leq 1,\ |z|\leq 1 $ and $ |w| \leq 1.$

According to Lemma \ref{1.7}, we may assume without any restriction that  $p_1\in [0,2]$. Thus , by applying the triangle inequality and taking $p_{1}= p$, we find that
\begin{equation}\label{eq3.19}
\begin{split}
&\left|a_{2}a_{4}-a_{3}^2\right| \leq \left|c_{2}c_{4}- c_{3}^2\right| + \frac{(1-\beta)}{2}\left|(c_{4}-c_{2}c_{3})\right|p \\
&\quad \quad\quad\quad +\frac{1}{4}(1-\beta) |c_{2}c_{3}|p + \frac{1}{8}(1-\beta)^2 \left|(c_{3}-10c_{2}^2)\right|p^2\\
&\quad \quad \quad\quad  +\frac{1}{2}|c_{2}|(1-\beta)^3 p^3 +\frac{1}{2}|c_{2}^3|(1-\beta)p+ \frac{(1-\beta)^4}{16}p^4\\
&\quad \quad \quad \quad  +\frac{(1-\beta)^2}{32}p^4 +\frac{1-\beta}{16}|c_{2}|p^3 + \frac{1}{4}|c_{2}^2|(1-\beta)p^2 +\frac{1}{8}|c_{2}|(1-\beta)^2 p^3\\
&\quad \quad \quad \quad + (|x|+|y|) \left[\frac{1-\beta}{24}|c_3- \frac{5}{4}c_{2}^2|(4-p^2)+ \frac{1}{24}|c_{2}|(1-\beta)^2 p(4-p^2) \right.\\
&\quad \quad \quad \quad + \frac{1}{96}(1-\beta)^3 p^2(4-p^2) + \frac{(1-\beta)^2}{32}p^2(4-p^2)+ \frac{(1-\beta)}{16}|c_{2}| p(4-p^2)\\
&\quad \quad \quad \quad + \left.\frac{1}{8}|c_{2}^2| (1-\beta)(4-p^2)+\frac{(1-\beta)}{16}c_{2}p (4-p^2)\right]\\
&\quad \quad \quad \quad +(|x|^2+|y|^2)\left[\frac{(1-\beta)^2}{64}p^2(4-p^2) +\frac{(1-\beta)}{32} |c_{2}|p(4-p^2)\right]\\
&\quad \quad \quad \quad + (1-|x|^2)\left[ \frac{1-\beta}{16}|c_{2}|(4-p^2)+ \frac{1-\beta)^2}{32}p(4-p^2)\right]\\
&\quad \quad \quad \quad + (1-|y|^2)\left[ \frac{1-\beta}{16}|c_{2}|(4-p^2)+ \frac{1-\beta)^2}{32}p(4-p^2)\right]\\
&\quad \quad \quad \quad + \frac{(1-\beta)^2}{144}(4-p^2)^2(|x|+|y|)^2.
\end{split}
\end{equation}
We now apply the Lemmas  \ref{1.6},  \ref{1.7} and  \ref{1.8}, and set $ \gamma_{1} =|x|\leq 1,\gamma_{2}=|y| \leq 1$. Then \eqref{eq3.19} can be rewritten in the form
$$|a_{2}a_{4}-a_{3}^2| \leq S_{1} + S_{2}(\gamma_{1} +\gamma_{2})+S_{3}(\gamma_{1}^2+\gamma_{2}^2)+S_{4}(\gamma_{1}+\gamma_{2})^2= F(\gamma_{1},\gamma_{2}), $$
where
\begin{eqnarray*}
S_{1} = S_{1}(p) &=& \frac{1}{8} + \frac{5}{6}(1-\beta)p+
\left[\frac{9}{8}(1-\beta)^2 + \frac{(1-\beta)}{4}\right]p^2 \\
& +& \left(\frac{1-\beta}{16}+\frac{(1-\beta)^3}{2}+\frac{(1-\beta)^2}{8}\right) p^3 \\
& +& \left[\frac{(1-\beta)^4}{16}+\frac{(1-\beta)^2}{32}\right]p^4 \\
& +&\frac{(1-\beta)}{8}(4-p^2)+\frac{(1-\beta)^2}{16}p(4-p^2) \geq 0, \\
S_{2} = S_{2}(p)&= & \left( \frac{(1-\beta)^2}{32} + \frac{(1-\beta)^3}{96}\right) p^2(4-p^2)\\
  &+&\left(\frac{1-\beta}{16}+\frac{5(1-\beta)^2}{48}\right)p(4-p^2)+\frac{(1-\beta)}{6}(4-p^2) \geq 0, \\
S_{3}= S_{3}(p) &= &\displaystyle \frac{(1-\beta)^2}{64}(p^2-4)(4-p^2) \leq 0, \\
S_{4}= S_{4}(p)&=& \displaystyle \frac{(1-\beta)^2}{144}(4-p^2)^2 \geq 0.
\end{eqnarray*}
Maximizing $F(\gamma_{1},\gamma_{2})$ in a square $\Delta:=\{(\gamma_{1},\gamma_{2}): 0 \leq \gamma_{1} \leq 1,0 \leq \gamma_{2} \leq 1\}$ we conclude  that $ \max F(\gamma_{1}, \gamma_{2})= F(1,1)$. Defining now a function $ K: [0,2] \rightarrow \mathbb{R}$  as in Theorem \ref{thm2.1} defined by
\begin{equation}\label{eq3.21}
K(p)= \max F(\gamma_{1},\gamma_{2})= F(1,1)= S_{1}+2S_{2}+2S_{3}+4S_{4},
\end{equation}
and analyzing its behavior, we infer that $K(p)$ is an increasing function of $p$ and attains its maximum  at $p=2$. Consequently, we have
$$\max_{0 \leq p \leq 2} K(p) = K(2) =\frac{1}{8}+ \frac{19}{6}(1-\beta)+6(1-\beta)^2+ 4(1-\beta)^3+(1-\beta)^4,$$
that completes the proof of the theorem.
\end{proof}

\begin{remark} For $\phi(z)=z$, Theorem \ref{thm3.1} reduces to the Theorem 1, due to \cite{Cag}.\end{remark}

\noindent\textbf{Conclusions}
 In the present paper, we have estimated  a smaller upper bound and more accurate estimation   for the functional $|a_2 a_4-a_3^{2}|$ for functions in the class of strongly  bi-close-to-convex functions of order $\alpha, (0 \leq \alpha \leq 1)$ and the class of bi-close-to convex functions of order $\beta, (0 \leq \beta <1)$. Obtaining a sharp estimate for $|a_2 a_4-a_3^{2}|$  in these classes are still open and keeps the researcher interested.\bigskip

\noindent\textbf{Acknowledgments}
The work of  the fourth author is  supported  by a grant from the Science and Engineering Research Board, Government of India under Mathematical Research Impact Centric Support of Department of Science and Technology (DST)(vide ref: MTR/2017/000607).\bigskip

\noindent\textbf{Availability of supporting data}

\noindent Not applicable.\bigskip

\noindent\textbf{Competing Interests}

\noindent The authors declare that they have no competing interests.\bigskip

\noindent\textbf{Authors' Contributions}

\noindent Each of the authors contributed to each part of this study equally, all authors read and approved the final manuscript.


\end{document}